\documentstyle[12pt]{article}
\title{\bf Exponential Stability of Linear Delay
Impulsive Differential Equations}
\author{A. Anokhin, \\ College of Physics and Mathematics, \\
452000, Chelyabinsk, Russia \\
L. Berezansky $^{1}$ \\
Department of Mathematics and Computer Science, \\
Beer-Sheva 84105, Israel, \\
E. Braverman $^{1}$ \\
Technion - Israel Institute of Technology,\\
Department of Mathematics, 32000, Haifa, Israel }

\begin{document}
\maketitle

\footnotetext[1]{Supported by : The Centre for
 Absorption in Science, Ministry of Immigrant Absorption State of Israel }

Corresponding author:
Elena Braverman
Technion - Israel Institute of Technology,
Department of Mathematics,
32000, Haifa, Israel

e-mail : maelena@tx.technion.ac.il

{\bf Abstract}

For an impulsive delay differential equation
with bounded delay and bounded coefficients
the following result is established.
If each solution is bounded on $[0, \infty )$
together with its derivative for each bounded
right-hand side then the equation is exponentially stable.
A coefficient stability theorem is presented.

\section{Introduction}

{}~~~~~The work of Millman and Myshkis [1] was the first one
dealing with impulsive differential equations.
Recently this field has been intensively investigated
(see the monographs [2,3] and their bibliography).
Not so much has been developed in the direction of impulsive
functional differential equations [4,5,6].

The paper deals with the exponential stability of a linear
impulsive delay differential equation. A new method
of research in stability is proposed.
The stability theorem is based on two results.

The first result is the representation of solutions.
The representation formulas were presented in our works
[7,8].
 However now we obtain the result that is more
convenient for an application.
The idea of its proof is also
different.
It is to be emphasized that a representation
of solutions is exploited for the stability investigation
not for a long time.
The idea was developed in the works of Corduneanu [9-11].

The second result is the Bohl-Perron theorem for an impulsive
equation.
The Bohl-Perron theorem is the result of the following
type.
Suppose that any solution of  the linear differential equation is bounded
on $[0, \infty )$ for each bounded on $[0, \infty)$
right-hand side.
Then the solution $X$ of the corresponding homogeneous equation
has an exponential estimate
$$ \parallel X(t) \parallel \leq N \exp (- \nu t), $$
with positive constants $N , \nu. $

For an ordinary differential equations without impulses
Bohl-Perron theorem is discussed in detail in [12].
For a delay differential equation results of this type
were obtained by Halanay [13] and Tyshkevich [14].
These results were generalized and completed in [15] and [16].
The exponential behavior of impulsive differential equations
was investigated in [17,18].

We prove the Bohl-Perron theorem for impulsive delay differential
equations on the base of the scheme proposed in [19].
The scheme was intensively used in the stability theory
of functional differential equations [20,21].

It is to be emphasized that in the impulsive conditions
$$ x(\tau _i + 0) = B_i x(\tau _i - 0) $$
we do not assume that $B_i$ are invertible.
{}From this point of view the result obtained is new
for impulsive equations without delay as well.
The point is that the equality for the Cauchy matrix
$$C(t,s) = X(t) X^{-1} (s) $$
holds for impulsive differential equations only if
$B_i$ are invertible, and the proof of the Bohl-Perron
theorem [17] is based on this equality.

In conclusion we present the exponential stability result.
On the base of this theorem sufficient stability conditions
for parameters of the equation are obtained.

\section{Preliminaries}

{}~~~~~Let $0 = \tau_0 < \tau_1 < \dots $ be the fixed points,
$\lim _{i \rightarrow \infty} \tau _i = \infty , $
${\bf R}^n$ be the space of $n$-dimensional column vectors
$x = col (x_1, \dots ,x_n)$ with the norm
$\parallel x \parallel = \max _{1 \leq i \leq n} \mid x_i \mid $,
by the same symbol $\parallel \cdot \parallel$ we shall
denote the corresponding matrix norm,

$E_n$ is an $n \times n$ unit matrix,

$\chi _e : [0, \infty ) \rightarrow {\bf R}
$ is the characteristic function of the set $e :
\chi_e (t) = 1$, if $t \in e, $ and $\chi_e (t) = 0, $
otherwise.

${\bf L}_{\infty} $ is a Banach space of essentially bounded
Lebesgue measurable functions $x: [0, \infty ) \rightarrow
{\bf R}^n , \parallel x \parallel _{{\bf L}_{\infty}} =
vraisup _ {t \geq 0} \parallel x(t) \parallel ,$

${\bf AC} $ is a linear space of functions $x: [0, \infty)
\rightarrow {\bf R}^n $ absolutely continuous on any segment
$[t, t+1],~ t \geq 0,$

${\bf PAC} (\tau _1, \dots , \tau_k, \dots ) $
is the space of piecewise absolutely continuous
functions $x: [0, \infty) \rightarrow {\bf R}^n,$
with jumps only at the points $\tau_1, \tau_2,
\dots, \tau _k, \dots, $  i.e.
$${\bf PAC}(\tau_1, \dots, \tau_k, \dots) =
\left\{ x:[0, \infty) \rightarrow {\bf R}^n \mid
x(t) = y(t) + \sum_{i=1}^{\infty} {\alpha _i
\chi _{[\tau _i, \infty )} (t)}, \right. $$ $$\left. t \geq 0,
\alpha_i \in {\bf R}^n, y \in {\bf AC} \right\}. $$

It is to be noted that a function $x \in {\bf PAC}
(\tau _1, \dots , \tau _k, \dots ) $ is right continuous.

We consider the linear delay differential equation
\begin{eqnarray}
\dot{x} (t) + \sum_{i=1}^m {A_i (t) x[h_i(t)]} =
r(t),~ t \geq 0,
\end{eqnarray}
\begin{eqnarray}
x(\tau _j) = B_j x(\tau _j - 0), j = 1,2, \dots,
\end{eqnarray}
\begin{eqnarray}
x(\xi ) = \varphi (\xi), \xi < 0.
\end{eqnarray}

Here $B_j$ are constant $n \times n$-matrices,
$h_i (t) \leq t, i = 1, \dots, m, t \geq 0.$

We consider the equation (1),(2),(3) under the
following assumptions

(a1) $0 = \tau_0 < \tau_1 < \tau _2 < \dots $ are fixed points,
$\lim_{i \rightarrow \infty} \tau _i  = \infty $ ;

(a2) $ r \in {\bf L}_{\infty} ,$ columns of $A_i , i = 1,
\dots , m,$ are in ${\bf L}_{\infty};$

(a3) $h_i : [0,\infty) \rightarrow {\bf R} $
are Lebesque measurable functions, $$h_i (t) \leq t, t \geq 0,
i = 1, \dots , m;$$

(a4) $\varphi : (- \infty, 0) \rightarrow {\bf R}^n $
is a Borel measurable bounded function;

(a5) $ M = \sup_i \parallel B_i \parallel
< \infty;$

(a6) $I =  \lim_{b \rightarrow \infty} \sup _{t,s >b}
\frac{i(t,s)}{t-s} < \infty. $

Here $i(t,s)$ is a number of points $\tau _j$ belonging
to the segment $[s,t].$
$$ $$

\underline{{\sl Definition 1.}}
 A function $x \in {\bf PAC}$
is said to be {\bf a solution of the impulsive equation
(1),(2)}  with the initial function $\varphi (t)$
if (1) is satisfied for almost all $t \in [0, \infty)$
and the equalities (2) hold.
$$  $$

Our objective is to study the exponential stability of
the impulsive equation (1),(2),(3), if any solution is bounded
for a bounded on $[0, \infty)$ right-hand side $r$.
$$ $$

\underline{ {\sl Definition 2.}}
The equation (1),(2),(3)
is said to be {\bf exponentially stable
} if there exist
positive constants $N$ and $\nu$ such that for any solution
$x$ of the corresponding homogeneous equation
\begin{eqnarray}
\dot{x} (t) + \sum_{i=1}^m {A_i (t) x[h_i(t)]} = 0
,~ t \geq 0,
\end{eqnarray}
(2),(3)
the inequality
\begin{eqnarray}
 \parallel x(t) \parallel \leq N \exp (- \nu t)
[\sup_{s<0} \parallel \varphi(s) \parallel
+ \parallel x(0) \parallel ]
\end{eqnarray}
holds.

If the inequality (5) is valid for $\varphi \equiv 0,$
then the equation is said to be {\bf exponentially stable
with respect to the initial value}, if (5) holds
for $x(0)=0,$ then the equation is said to be {\bf exponentially stable
with respect to the initial function}.
$$ $$

It is to be emphasized that exponential properties
of delay impulsive differential equations differ
greatly from the properties of the corresponding
equations without delay (or without impulses).
For instance in [22] the following theorem
was proven in the scalar case.

Suppose that the Cauchy problem $\dot{x} + Ax = f, x(0)=0$
has a bounded on $[0, \infty ) $ solution for any $f \in {\bf
 L}_{\infty} $ and there exists a positive constant $M_1$
such that $$\mid B_p B_{p+1} \dots B_j \mid \leq M_1 $$
for any positive integers $p,j, p<j, A \in {\bf L}_{\infty}.$
Then there exist positive constants $N$ and $\nu$ such
that the solution of the corresponding homogeneous equation
satisfying (4),(2), $x(0)=1,$ has the exponential estimate
$$ \mid x(t) \mid \leq N \exp (- \nu t). $$

The following example shows that for the delay impulsive
differential equation this theorem is not valid.

\underline{{\sl Example.}} The equation
$$ \dot{x} (t) + a x(t-1) = f, ~~ 0<a<\pi /2,$$
is exponentially stable
[23].
Consider this equation with impulsive conditions
$$x(j) = - x(j-1),~ j = 1,2, \dots .$$
Here $ \mid B_j \mid = 1.$ The  absolute
value of the solution of the corresponding homogeneous equation with
$x(0)=1$ is steadily growing, $t >1$
$$ \mid x(t) \mid = 1,~ t \in [0,1), $$
$$ \mid x(t) \mid \geq 1 + a(t - 1) ,~ t \in [1,2), $$
$$ \mid x(t) \mid \geq 1 + a +
a(1+a)(t-2),~ t \in [2,3),$$
$$ \mid x(t) \mid \geq (1+a)^2 + a (1+a)^2 (t-3), t \in [3,4), $$
$\dots .$

\section{Representation of solutions}

{}~~~~~The main result of this section deals with the representation
of solutions.

In the stability theory of ordinary differential equations the
representation of the Cauchy matrix
$$ C(t,s) = X(t) X^{-1} (s) $$
is intensively used.
Here $C$ is the Cauchy matrix, i.e. the kernel of the integral
representation of solutions, $X(t)$ is the solution of the corresponding
homogeneous equation satisfying $X(0) = E_n.$

For delay differential equations this equality
generally speaking is not true.
$$ $$

\underline{\sl Definition 3.}
 An impulsive delay differential
equation
\begin{eqnarray}
\dot{x} (t) + \sum_{i=1}^m {A_i (t) x[h_i(t)]} =
0 ,~ t \geq s,
\end{eqnarray}
\begin{eqnarray}
x(\xi) = 0, ~ \xi < s,
\end{eqnarray}
$$x(\tau _j) = B_j x(\tau _j - 0), ~ j = 1,2, \dots, \tau _i > s, $$
is said to be {\bf a homogeneous "${\bf s}$ - curtailed" equation }
(1),(2).
The solution $X(\cdot,s)$ of this equation satisfying $X(s,s) = E_n$
is said to be {\bf the fundamental matrix of
the "$s$ - curtailed" equation.}

The main result of this section is the following.

\newtheorem{guess}{Theorem}[section]
\begin{guess}
Suppose (a1)-(a4) are satisfied.
Then the solution of the problem (1),(3) , with the initial
value
\begin{eqnarray}
x(0) = \alpha _0
\end{eqnarray}
and impulsive conditions
\begin{eqnarray}
x(\tau _i) = B_i x(\tau _i - 0) + \alpha _i ,~ i = 1,2, \dots,
\end{eqnarray}
can be presented as
\begin{eqnarray}
x(t) = \int_0^t {X(t,s) r(s) ds}
- \sum_{i=1}^m {\int_0^t {X(t,s) A_i(s) \varphi (h_i(s)) ds}} +
\sum_{j=0}^{\infty} {X(t, \tau_j) \alpha _j} .
\end{eqnarray}
Here $\varphi (\zeta) = 0,$ if $\zeta \geq 0. $
\end{guess}

First we prove that there exists one and only one
solution of the initial problem for the equation (1),(9),(3)
(Lemma 3.1).
Then we establish the coincidence of the Cauchy matrix and
the fundamental matrix of "$s$-curtailed" equation (Lemma 3.3).
To this end we need the estimate of the latter matrix (Lemma 3.2).
Afterwards we prove the Theorem 3.1.

\newtheorem{uess}{Lemma}[section]
\begin{uess}
Suppose that the hypotheses
(a1) - (a4) hold.
 Then the Cauchy problem for the delay impulsive
differential equation (1),(3),(8),(9) has
a unique solution in ${\bf PAC} (\tau _1, \dots, \tau _k , \dots )$ .
\end{uess}

{\sl Proof.} Under the hypotheses (a1) -(a4) there exists
one and only one solution of the delay differential equation (1),(3)
on the interval $[0, \tau _1 )$ [23].
 The solution of the impulsive
equation (1),(3) on the interval $[\tau _1, \tau _2)$
can be treated as the Cauchy problem for the delay
differential equation (without impulses) on this interval
with the initial value
$$ x(\tau _1) = B_1 x(\tau _1 - 0) + \alpha_1  $$
and the initial function
$$ \varphi _1 (s) = \varphi (s), s < 0,
\varphi _1 (s) = x(s), s \in [0, \tau _1). $$
This initial function is bounded and Borel measurable
since $\varphi$ possesses this property and $x$ is
continuous and bounded on $[0, \tau _1).$
By induction the solution can be constructed on
the whole semi-axis $[0,\infty)$ and it is obviously unique.
The proof of the lemma is complete.
$$ $$

In the stability theory for delay differential equations
the "$s$ -curtailment" theorem is used:
if $C(t,s)$ is the Cauchy matrix of the equation (1),
then $C(\cdot ,s)$ is the solution of the Cauchy problem
(6),(7), $x(s) = E_n.$

For impulsive delay differential equations the same statement is valid.
First an auxiliary assertion will be proven.

\begin{uess}
 For the fundamental matrix of
"$s$ -curtailed" equation the following estimate is valid
$$\parallel X(t,s) \parallel \leq  \prod _{s< \tau_i \leq t}
(1 + \parallel B_i \parallel ) \exp \left\{ \int_s^t { \sum _{k=1}^m
{ \parallel A_k (\zeta) \parallel d \zeta }} \right\}. $$
\end{uess}

{\sl Proof.}
 Let $s < \tau _i < \tau _{i+1} < \dots < \tau _j \leq t.$
Then for $t\in [s, \tau _i) $ the solution $x$ of the problem
(6),(7),(2), $x(s) = E_n$ can be presented as
$$ x(t) = E_n - \int_s^t { \sum _{k=1}^m {A_k (\zeta)
x(h_k (\zeta)) d \zeta }} . $$
Hence
\begin{eqnarray}
 \parallel x(t) \parallel \leq
1 + \int_s^t { \sum _{k=1}^m {\parallel A_k (\zeta) \parallel
\sup _{\xi \in [s,\zeta ]} { \parallel x(\xi ) \parallel } d \zeta }} .
\end{eqnarray}
Denote $y(t) =
\sup _{\zeta \in [s,t]} { \parallel x(\zeta) \parallel }  . $
For the function $y(t)$ the inequality (11) implies
$$y(t) \leq 1 + \int_s^t { \sum _{k=1}^m {\parallel A_k (\zeta) \parallel
y(\zeta)  d \zeta }} .$$
Then the Gronwall - Bellman inequality gives
$$y(t) \leq \exp
\left\{ \int_s^t { \sum _{k=1}^m {\parallel A_k (\zeta) \parallel
d \zeta }} \right\} .$$
Therefore for the solution $x$ of the problem (6),(7),(2), $x(s) = E_n$
we have obtained the estimate
\begin{eqnarray}
\parallel x(t) \parallel  \leq \exp \left\{ \int_s^t { \sum _{k=1}^m
 {\parallel A_k (\zeta) \parallel
d \zeta }} \right\} .
\end{eqnarray}
Let $\tau _i \leq t < \tau _{i+1} .$
Then
$$ x(t) = x(\tau _i) -  \int_{\tau _i}^t { \sum _{k=1}^m {A_k (\zeta)
x(h_k (\zeta)) d \zeta }} . $$
Thus the inequality (12) and the impulsive condition
$x(\tau _i) = B_i x(\tau _i - 0)$ imply the estimate
$$\parallel x(t) \parallel  \leq  ( 1 + \parallel B_i \parallel )
\exp \left\{ \int_s^{\tau _i} { \sum _{k=1}^m
 {\parallel A_k (\zeta) \parallel
d \zeta }} \right\} + $$
$$ + \int_{\tau _i}^t { \sum _{k=1}^m {\parallel A_k (\zeta) \parallel
 { \parallel x[h_k (\zeta)] \parallel } d \zeta }} .
$$
Hence
if we again denote $ y(t) = \max _{\zeta \in [s,t]}
\parallel x(\zeta) \parallel $ we obtain
$$ y(t)  \leq  ( 1 + \parallel B_i \parallel )
\exp \left\{ \int_s^{\tau _i} { \sum _{k=1}^m
 {\parallel A_k (\zeta) \parallel
d \zeta }} \right\} + $$
$$ + \int_{\tau _i}^t { \sum _{k=1}^m {\parallel A_k (\zeta) \parallel
 { y(\zeta) } d \zeta }} .
$$

Repeating the previous argument gives
$$\parallel x(t) \parallel  \leq (1 + \parallel B_i \parallel)
\exp \left\{ \int_s^{\tau _i} { \sum _{k=1}^m
 {\parallel A_k (\zeta) \parallel
d \zeta }} \right\}
\exp \left\{ \int_{\tau _i}^t { \sum _{k=1}^m {\parallel A_k (\zeta) \parallel
 d \zeta }} \right\} = $$
$$ = ( 1 + \parallel B_i \parallel )
\exp \left\{ \int_s^t { \sum _{k=1}^m
 {\parallel A_k (\zeta) \parallel
d \zeta }} \right\} , ~ t \in [\tau _i, \tau _{i+1} ) . $$
By considering the solution $x$ of the problem (6),(7), $x(s) = E_n$
in the intervals $[\tau _{i+1}, \tau _{i+2} ), \dots ,
[\tau _{j-1} , \tau _j )$ and at the point $\tau _j$
we obtain the required inequality for $\parallel X(t,s) \parallel $.
This completes the proof of the lemma.
$$ $$

\underline{\sl Remark.}
If $B_i \neq 0, i= 1,2, \dots, $ then in the statement of
the lemma one can write
$\prod \parallel B_i \parallel $ instead of
$\prod (1 + \parallel B_i \parallel). $
$$ $$

\underline{\sl Corollary.}
For any $b>0$ the function $X(t,s)$ is bounded in $[0,b] \times [0,b]$.
$$ $$

This result immediately follows from the inequality
$$\parallel X(t,s) \parallel \leq  \prod _{0 < \tau_i \leq b}
(\parallel B_i \parallel + 1 ) \exp \left\{ \int_0^b { \sum _{k=1}^m
{ \parallel A_k (\zeta) \parallel d \zeta }} \right\}, $$
where $t \leq b.$
$$ $$

Now we can prove the "$s$ - curtailment" result for impulsive
delay differential equations.

\begin{uess}
Let $X(t,s)$ be the fundamental matrix of
the "$s$ -curtailed" equation . Then the solution $y$ of the Cauchy
problem (1),(2) , $ y(\xi) = 0, \xi~<~0, y(0) = 0 $
can be presented as
\begin{eqnarray}
y(t) = \int_0^t {X(t,s) r(s) ds}.
\end{eqnarray}
\end{uess}

{\sl Proof.}
We shall prove that (13) is the solution of the
problem (1),(2) , $ y(\xi) = 0, \xi < 0, y(0) = 0 $.
By differentiating the equality (13) in $t$ we obtain
\begin{eqnarray}
\dot{y}(t) = r(t) + \int_0^t {X'_t (t,s) r(s) ds},
\end{eqnarray}
since $X(t,t) = E_n .$
The equality (13) implies
$$y[h_k (t)] = \int_0^{h_k^+ (t)} {X(h_k (t),s) r(s) ds} = $$
$$= \int_0^t {X(h_k (t),s) r(s) ds} -
\int_{h_k^+ (t)}^t {X(h_k (t),s) r(s) ds} , $$
where $h^+ = max\{h, 0\} .$

As $X(t,s) = 0 $ for $t < s $ then the second integral
in the right-hand side vanishes.
Hence
$$ y(h_k (t)) = \int_0^t {X(h_k (t),s) r(s) ds} ,$$
that together with the equality (14) gives that $y$
is the solution of the problem
(1) , $ y(\xi) = 0,~ \xi < 0,~ y(0) = 0 $.

It remains to show that $y$ satisfies the
impulsive conditions
$$y(\tau_i) = B_i y(\tau _i - 0). $$
Let $i$ be a fixed positive integer and
$\{t_k \}_{k=1}^{\infty} \subset [0, \tau _i) $
be such sequence that $t_k$ tends to $\tau _i$ as
$k \rightarrow \infty .$
We shall prove that the equality
\begin{eqnarray}
\lim _{t_k \rightarrow \tau _i - 0} \int_0^{t_k} {X(t_k,s) r(s) ds}
= \int_0^{\tau _i} {X(\tau _i - 0,s) r(s) ds}
\end{eqnarray}
holds, i.e that the limit under the integral is possible.

Denote
$$g_k(s) = X(t_k,s) \chi _{[0,t_k]} (s) r(s), $$
$$g(s) = X(\tau _i - 0, s) \chi _{[0, \tau _i) } (s) r(s).$$
Evidently
$$X(t_k,s) \rightarrow X(\tau _i - 0,s)$$
and
$$\mid \chi _{[0,t_k)} (s) - \chi _{[0,\tau_i)} (s) \mid =
\chi _{(t_k, \tau _i)} \rightarrow 0 $$
as $k \rightarrow \infty$ for any $s \geq 0.$

Therefore
$\lim _{k \rightarrow \infty} {g_k (s)} = g(s), s \geq 0.$
Besides this, Lemma 3.2 gives (see the corollary)
$$\parallel g_k (s) \parallel  \leq
\prod _{0 < j < i}
(\parallel B_j \parallel + 1 ) \exp \left\{ \int_0^{\tau _i} { \sum _{k=1}^m
{ \parallel A_k (\zeta) \parallel d \zeta }} \right\}
\parallel r(s) \parallel . $$
Therefore the functions $g_k (s)$ are uniformly bounded
for $s \leq \tau _i .$

By the Lebesgue theorem on limit under the integral
we obtain (15).
The function $X(t,s)$ satisfies the impulsive condition
$X(\tau _i, s) = B_i X(\tau _i - 0, s) .$
Thus the equality (15) implies
$$ B_i y(\tau _i - 0)
=     B_i  \lim _{t_k \rightarrow \tau _i - 0}
 \int_0^{t_k} {X(t_k,s) r(s) ds} =
$$ $$= \int_0^{\tau _i} {B_i X(\tau _i - 0,s) r(s) ds}
=
\int_0^{\tau _i} {X(\tau _i ,s) r(s) ds}
= y(\tau _i) .
$$
Hence
$ y(\tau _i) =
B_i y(\tau _i - 0) , $
which completes the proof of the lemma.
$$ $$

\underline{\sl Proof of Theorem 3.1.}
By Lemma 3.1
there exists a unique solution of the Cauchy problem
(
we notice that the sum
$\sum_{j=0}^{\infty} {X(t, \tau_j) \alpha _j} $
is definite since for each $t > 0 $ this sum contains only a finite
number of terms with $\tau _j \leq t).$
By the direct substitution one can be convinced
that the solution of the problem
\begin{eqnarray}
\dot{x} (t) + \sum_{k=1}^m {A_k (t) x[h_k (t)] } =
r(t) - \sum_{k=1}^m {A_k (t) \varphi [h_k (t)] }, \nonumber \\
x(\xi) =  0, ~ \xi < 0, ~
\varphi (\zeta) = 0, ~ \zeta \geq 0,
\end{eqnarray}
 coincides with the solution of the problem
(1), $x(0)=0.$

Lemma 3.3 gives that the solution of the
problem (16),(2), $x(\xi) = 0, \xi \leq 0 $
,$x(0)=0$ can be presented as
$$ x_1 (t) = \int_0^t {X(t,s) r(s) ds}
- \sum_{i=1}^m {\int_0^t {X(t,s) A_i(s) \varphi (h_i(s)) ds}} .$$

Define $x_2 (t) =
\sum_{j=0}^{\infty} {X(t, \tau_j) \alpha _j} .$

Since $X(t,s)$ is the fundamental matrix of the "$s$-
curtailed" equation then
$$X'_t (t, \tau_j) \alpha _j +
\sum_{i=1}^m {A_i (t) X(t, h_i( \tau_j))
 \alpha _j} = 0 ,~ j = 1,2, \dots .$$
Thus
$$
\sum_{j=0}^{\infty} {
X'_t (t, \tau_j) \alpha _j} +
\sum_{j=0}^{\infty}{
\sum_{i=1}^m {A_i (t) X(t, h_i (\tau_j)) \alpha _j}} = 0
 .$$
Therefore $x = x_1 + x_2 $
satisfies the equation (16).

It remains to show that the initial condition
and the impulsive conditions (9) are satisfied.

For any $s \geq 0$, $j = 1,2, \dots $
we have (see the proof of Lemma 3.3)
$$ X(\tau _j, s) = B_j X(\tau_j - 0, s)$$
and
$$ \int_0^{\tau _j} {X(\tau _j, s) g(s) ds}
= B_j \int_0^{\tau _j} {X(\tau_j - 0, s) g(s) ds} $$
for any $g \in {\bf L}_{\infty} .$

Hence
$$x(\tau _j) = \int_0^{\tau _j} {X(\tau _j,s) r(s) ds} - $$ $$
- \sum_{i=1}^m {\int_0^{\tau _j} {X(\tau _j,s) A_i(s) \varphi (h_i(s)) ds}} +
\sum_{k=0}^{j} {X(\tau _j, \tau_k) \alpha _k}=
$$  $$
 = B_j \int_0^{\tau _j} {X(\tau _j - 0,s) r(s) ds} -
 \sum_{i=1}^m {B_j \int_0^{\tau _j}
 {X(\tau _j  - 0 ,s) A_i(s) \varphi (h_i(s)) ds}} + $$ $$
+ B_j \sum_{k=0}^{j - 1} {X(\tau _j - 0 , \tau_k) \alpha _k} +
X(\tau _j, \tau _j) \alpha _j = $$ $$
 = B_j [ \int_0^{\tau _j} {X(\tau _j - 0,s) r(s) ds} -
 \sum_{i=1}^m { \int_0^{\tau _j}
 {X(\tau _j  - 0 ,s) A_i(s) \varphi (h_i(s)) ds}} + $$ $$
+  \sum_{k=0}^{j - 1} {X(\tau _j - 0 , \tau_k) \alpha _k}] +
E_n \alpha _j = B_j x(\tau _j - 0) + \alpha _j, $$
i.e the impulsive conditions (9) are satisfied.

Further,
$$x(0) = \sum_{i=0}^{\infty} {X(0, \tau _i) \alpha _i} =
X(0,0) \alpha _0 = \alpha _0 $$
since $X(t,s) = 0 , ~ t < s.$
Therefore the initial condition (8) is also satisfied.
The proof of the theorem is complete.
$$ $$

\underline{\sl Definition 4.}
The operator $C: {\bf L}_{\infty}
\rightarrow {\bf PAC} (\tau_1 , \dots , \tau_k , \dots ) $
defined by the formula
$$ (Cf)(t) = \int_0^t {X(t,s) f(s) ds} $$
is said to be {\bf the Cauchy operator }
of the impulsive delay differential equation (1),(2),(3).

\section{Exponential estimates of the Cauchy matrix}

{}~~~~~The purpose of this work is to obtain
the exponential estimate of the matrix $X(t,s)$
and to investigate the exponential stability
of the impulsive delay differential equations.

Let ${\bf D}_{\infty}  \subset {\bf PAC} (\tau_1 , \dots,
\tau_k , \dots ) $ be a space of functions
$x:[0, \infty )\rightarrow {\bf R}^n$
 absolutely continuous on the intervals $[\tau_j , \tau _{j+1} )$
satisfying the impulsive conditions  (2)
and such that both $x$ and its derivative are essentially
bounded on $[0, \infty ). $

We introduce the norm in ${\bf D}_{\infty}$
$$ \parallel x \parallel _{{\bf D}_{\infty}} =
\parallel x \parallel _{{\bf L}_{\infty}}
+ \parallel \dot{x} \parallel _{{\bf L}_{\infty}}. $$

The main result of this section is the following.
$$ $$

\begin{guess}
Let exist constant $\delta > 0$  such that
$$\theta _i (t) = t - h_i (t) \leq \delta, t \geq 0, $$
and the hypotheses (a1) -- (a6) hold.
Suppose that for any essentially bounded on $[0, \infty )$
right-hand side $r $ all solutions of the impulsive equation
(1),(2),(3) with
$\varphi \equiv 0,$
 are essentially bounded on $[0, \infty)$ together with their derivatives.

Then there exist positive constants $N$
and $\nu$ such that for the Cauchy matrix of the impulsive equation
(1),(2),(3) the inequality
$$ \parallel X(t,s) \parallel \leq N \exp [ - \nu (t-s)]$$
holds.
\end{guess}
$$ $$

For proving this theorem we need some auxiliary assertions.
$$ $$

\begin{uess}
 Suppose the hypotheses
(a5) and (a6) hold, i.e.
\\$ M = \sup_i \parallel B_i \parallel
< \infty$  and
$I =  \lim_ {b \rightarrow \infty} \sup _{t,s > b}
\frac{i(t,s)}{t-s} < \infty, $ \\
where  $i(t,s)$ is a number of the points $\tau_j$
belonging to the segment  $[s,t]$.
{\it Then } ${\bf D}_{\infty}$ {\it is a Banach space.}
\end{uess}

{\sl Proof.} We choose a positive constant $a$ such that
$I \ln M < a.$
Then for the linear impulsive equation
$\dot{x} + ax = f,$ (2)
the Cauchy matrix $C_0 (t,s)$
has an exponential estimate.
For proving this we use the representation [2,3]
$$C_0 (t,s) = \exp [- a(t-s)] \prod _{s<\tau_i \leq t} {B_i}.$$
Thus the following estimate is valid
$$\parallel C_0 (t,s) \parallel \leq
\exp \left\{- \left[a - \frac{\ln M \cdot i(t,s)}{t-s}
\right] (t-s) \right\} \leq
\exp [-\nu (t-s)], $$
where $\nu = a - I \ln M.$

Denote
$$({\cal L} x)(t) = \dot{x} + ax, $$
with $x(0) = 0,~ x(\tau_i) = B_i x(\tau_i - 0),$
$$(Cf)(t) = \int_0^t {C_0(t,s) f(s)ds}.
$$

Thus we have introduced the linear operators
${\cal L} : {\bf D}_{\infty} \rightarrow {\bf L}_{\infty}$
and $C : {\bf L}_{\infty} \rightarrow {\bf D}_{\infty}.$
Now we shall prove that these operators are continuous.

To this end
$$ \parallel {\cal L} x \parallel _{{\bf L}_{\infty}}=
 \parallel \dot{x} + ax \parallel _{{\bf L}_{\infty}}
 \leq  (1 + a) \parallel x \parallel _{{\bf D}_{\infty}},$$

$$ \parallel Cf \parallel _{{\bf D}_{\infty}}=
 \parallel x \parallel _{{\bf L}_{\infty}}
 + \parallel \dot{x} \parallel _{{\bf L}_{\infty}},
$$
where
 $$\parallel x \parallel _{{\bf L}_{\infty}}
 \leq \frac{1}{\nu} \parallel f \parallel _{{\bf L}_{\infty}},$$
 $$ \parallel \dot{x} \parallel _{{\bf L}_{\infty}} =
 \parallel f - a x \parallel _{{\bf L}_{\infty}} \leq $$
 $$\leq \parallel f \parallel _{{\bf L}_{\infty}}
 + a \frac{1}{\nu} \parallel f \parallel _{{\bf L}_{\infty}} =
 \left(1 + \frac{a}{\nu} \right) \parallel f \parallel _{{\bf L}_{\infty}}, $$
therefore
$$ \parallel Cf \parallel _{{\bf D}_{\infty}} \leq
 (1 + \frac{1 + a}{\nu}) \parallel f \parallel _{{\bf L}_{\infty}}. $$

Now we prove that the space ${\bf D}_{\infty}$ is complete.
Let $\{ x_k \}_{k=1}^{\infty} \subset {\bf D}_p $
be a sequence such that
$ \parallel x_k - x_i \parallel _{{\bf D}_{\infty}}
\rightarrow 0 $
as $k,i \rightarrow \infty.$
 Consider the functions
$f_k = {\cal L} x_k$.
Then
$$x_k (t)  = C_0 (t,0) +
\int_0^t {C_0 (t,s) f_k (s) ds} , $$
where $C_0 (t,0) $ is the solution of the
problem
$$\dot{x} + ax = 0, ~ x(0) = E_n , ~x(\tau_i) =
B_i x(\tau_i - 0) .$$

First we show that
the sequence $\{ g_k (t) \}$ , with
$g_k (t) = C_0(t,0)x_k(0), $
converges in ${\bf D}_{\infty} $.
Let $b < \tau_1 $.
Then
$ \parallel x_k - x_i \parallel _{{\bf D}_{\infty}}
\rightarrow 0 $  as $k,i \rightarrow \infty $ implies
$$ \max _{t \in [0,b]} [
\parallel x_k (t) - x_i (t) \parallel ]
 + vraisup _{t \in [0,b]} [
\parallel \dot{x}_k (t) - \dot{x}_i (t) \parallel ]
\rightarrow 0 . $$
Therefore
$$ \parallel x_k (0) - x_i (0) \parallel
\rightarrow 0 . $$
Since ${\bf R}^n $ is complete then there exists
$\beta$ such that
$x_k (0) \rightarrow \beta $ as $k \rightarrow \infty.$

The estimates $\parallel C_0(t,0) \parallel
\leq \exp (-\nu t) $   and
$\parallel \dot{C}_0(t,0) \parallel
\leq a \exp (-\nu t) $
imply that
$g_k(t) \rightarrow C_0(t,0) \beta $ in
${\bf D}_{\infty}. $

Denote
$$ \tilde{x}_k = \int_0^t {C_0(t,s) f_k (s) ds} . $$
Then
$$ \parallel \tilde{x}_k - \tilde{x}_i \parallel _{{\bf D}_{\infty}}
= \parallel \int_0^t {C_0(t,s) [f_k (s) - f_i (s)] ds} \parallel _{{\bf D}
_{\infty}}
\rightarrow 0 $$
 implies
$$ \parallel f_k  - f_i  \parallel _{{\bf L}
_{\infty}}
= \parallel {\cal L} (x_k - x_i) \parallel _{{\bf L}_{\infty}}
\rightarrow 0 $$
since ${\cal L} : {\bf D}_{\infty} \rightarrow
{\bf L}_{\infty}$
is a continuous operator.

As ${\bf L}_{\infty}$ is complete then there exists
$f \in {\bf L}_{\infty} $
such that
$$ \lim_{k \rightarrow \infty} { \parallel f_k  - f  \parallel _{{\bf L}
_{\infty}}} = 0. $$

Let $\tilde{x} = Cf.$ Since $C$ is a continuous operator, then
$$ \parallel \tilde{x}_k - \tilde{x} \parallel _{{\bf D}_{\infty} }
= \parallel C (f_k - f )  \parallel _{{\bf D}
_{\infty} }
\rightarrow 0 . $$
Thus $\tilde{x} = \lim _{k \rightarrow \infty} {\tilde{x}_k }
\in {\bf D}_{\infty}.$

Let $$x = C_0  (t,0) \beta + \tilde{x} . $$
Then $x_k$ tends to $x$ in ${\bf D}_{\infty}$,
which completes the proof of the lemma.
$$ $$

\begin{uess}
Let ${\cal L} , {\cal M} :
{\bf D}_{\infty} \rightarrow
{\bf L}_{\infty} $ be linear bounded operators
and $C , C_{\cal M}$ be the Cauchy operator of the equations
${\cal L} x = f , {\cal M } x = f $ correspondingly.
 Suppose  the Cauchy operator $C$ is a bounded operator acting from
$ {\bf L}_{\infty}$ to $  {\bf D}_{\infty}$ and the operator $
{\cal M} C : {\bf L}_{\infty} \rightarrow {\bf L}_{\infty}$
is invertible.
Then the Cauchy operator $C_{\cal M}$
also
acts from ${\bf L}_{\infty}$ to ${\bf D}_{\infty}$
and is bounded.
\end{uess}

{\sl Proof.}
For each $ f \in L_{\infty} $ the function
$x = C ({\cal M} C)^{-1} f $
is the solution of the problem $Mx = f, x(0) = 0.$
By the definition of the Cauchy operator
$C_{{\cal M}} = C ({\cal M} C )^{-1}.$
This implies that $C_{\cal M}$ is a bounded operator
acting from
${\bf L}_{\infty}$ to ${\bf D}_{\infty}$,
which completes the proof of the lemma.
$$ $$

\underline{\sl Proof of Theorem 4.1.}
First we shall obtain an exponential estimate
for $X(t,0).$
Let $\nu$ be a positive number.
Denote
$$y(t) = x(t) \exp {(\nu t )}.$$
Everywhere below we assume $x(\xi) = y(\xi) = 0, \xi < 0.$
If $x(t)$ satisfies the impulsive conditions (2), then
evidently
$$ y(\tau _j) = B_j y(\tau _j - 0), ~ j = 1,2, \dots. $$
We denote
$$ ({\cal L} x)(t) =
\dot{x} (t) + \sum_{i=1}^m {A_i (t) x[h_i (t)]}. $$

By substituting $x(t) = y(t) \exp {(-\nu t)} $ we obtain
$$ ({\cal L} x)(t) =
\exp {(- \nu t)} \dot{y} (t) -
\nu \exp {(- \nu t)} y(t) + $$ $$
+ \sum_{i=1}^m {
\exp {( - \nu h_i (t) )}
A_i (t) y[h_i (t)]} =
\exp {(- \nu t)} \left\{ \dot{y} (t) +
\sum_{i=1}^m {A_i (t) y[h_i (t)]} + \right. $$
$$
+ \left. \sum_{i=1}^m {
\exp {[  \nu (t - h_i (t) )] }
A_i (t) y[h_i (t)]} -
\sum_{i=1}^m {A_i (t) y[h_i (t)]} - \nu y(t) \right\} = $$
$$ = \exp {(- \nu t)} \left\{ ({\cal L} y) (t) - \nu y(t) +
\sum_{i=1}^m {
[ \exp (  \nu \theta _i (t))  - 1 ]
A_i (t) y[h_i (t)]}  \right\}.  $$
Denote
$$ ({\cal T} y)(t) =
\sum_{i=1}^m {[
\exp {(  \nu \theta _i (t) )} - 1 ]
A_i (t) y[h_i (t)]}  -  \nu y(t),   $$
$$({\cal M} y)(t) = ({\cal L} y)(t) + ({\cal T} y)(t).$$
Then
$$({\cal L} x)(t) = \exp {(-\nu t)} ({\cal M} y)(t) .$$

The impulsive equation (1),(2),(3)
 has a solution that can be presented as (13).
Let $C$ be the Cauchy operator of the equation (1),(2).

Under the hypotheses of the theorem the operator
${\cal L} : {\bf D}_{\infty} \rightarrow {\bf L}_{\infty}$
is bounded.
Since a solution  is in ${\bf L}_{\infty}$
together with its derivative for any
right-hand side from ${\bf L}_{\infty}$
then the Cauchy operator $C$ acts from
${\bf L}_{\infty}$ to
${\bf D}_{\infty} .$

Denote ${\bf D}_{\infty}^0 =
\{ x \in {\bf D}_{\infty} : x(0)= 0\}.$
Let ${\cal L}^0 $ be the contraction of
${\cal L}$ to ${\bf D}_{\infty}^0.$
Then the operator ${\cal L}^0~:{\bf D}_{\infty}^0
\rightarrow {\bf L}_{\infty} $
acts onto the space ${\bf L}_{\infty}.$
Therefore by the Banach theorem the Cauchy operator
$C : {\bf L}_{\infty} \rightarrow
{\bf D}_{\infty}^0 $ is bounded as its inverse.

By (a2) the columns of  $A_i$ belong
to $ {\bf L}_{\infty}$, i.e. there exists
$Q>0$ such that
$$ vraisup _{t \geq 0} \sum_{i=1}^m {\parallel A_i (t) \parallel} <Q.$$
Therefore the operator ${\cal T} : {\bf D}{\infty}_
\rightarrow {\bf L}_{\infty} $ is bounded, with
$$
\parallel {\cal T} \parallel _{
 {\bf D}{\infty}_
\rightarrow {\bf L}_{\infty}} \leq Q [ \exp {(\nu
\delta)} - 1 ] + \nu . $$
Thus ${\cal M} C = {\cal L} C + {\cal T}C =
E + {\cal T} C $ has a bounded inverse operator whenever
$$\parallel {\cal T} C \parallel _{
 {\bf L}{\infty}_
\rightarrow {\bf L}_{\infty}} < 1. $$
Here $E$ is the identity operator.

Let
$P =
\parallel C \parallel _{
 {\bf L}{\infty}_
\rightarrow {\bf D}_{\infty}} . $
Then
$$\parallel {\cal T} C \parallel _{
 {\bf L}{\infty}_
\rightarrow {\bf L}_{\infty}    } \leq P
\parallel {\cal T}  \parallel _{
 {\bf D}_{\infty}
\rightarrow {\bf L}_{\infty}     } \leq P
Q [\exp {(\nu \delta)} - 1] + P \nu < 1 $$
for $\nu$ being small enough.

The operators
${\cal L}$
and ${\cal T}$
continuously act from
${\bf D}_{\infty}$  to ${\bf L}_{\infty}.$
Therefore the operator
${\cal M} = {\cal L} - {\cal T}$
also continuously acts from
${\bf D}_{\infty}$ to ${\bf L}_{\infty}.$
Hence Lemma 4.2 implies that the Cauchy operator
$C_{\cal M} $ of the equation ${\cal M} y = f$, $x(0)=0$
continuously acts from
${\bf L}_{\infty}$ into ${\bf D}_{\infty}.$

We introduce
$$
\Psi _0 (t) =
 \left\{ \begin{array}{ll} E_n , & \mbox{if}~ t \leq \tau _1 , \\
\prod _{0 <  \tau _j \leq t} {B_j}, &
\mbox{if}~ t> \tau _1,~ M \leq 1, \\
\exp \{ - \ln {M} (I + \varepsilon ) t\}
\prod _{0 <  \tau _j \leq t} {B_j}, &
\mbox{if}~ t > \tau_1,~ M>1,
\end{array}
\right.  $$
where $\varepsilon$ is a positive constant,
$I$ and $M$ are the numbers defined in
the hypotheses (a5),(a6).

Let us prove that columns of both $\Psi _0 (t)$
and its derivative  are in ${\bf L}_{\infty}$.

Let $M \leq 1$.
Evidently
$\parallel \Psi _0 (t) \parallel \leq 1 $
and
$\dot{\Psi}_0(t)=0 , $
 therefore in the case $M \leq 1$
the columns of  $\Psi $ and $\dot{\Psi} $ are in
${\bf L}_{\infty} . $

Let $M>1$.
By the definition of the function $i(t,s)$
for any  fixed $\varepsilon > 0$ there exists
$b>0$ such that
$$ I > \frac{i(t,s)}{t-s} - \varepsilon $$
holds for any $t,s \geq b.$

Therefore for $t>b$
$$ \parallel \Psi _0 (t) \parallel
\leq \prod _{0 < \tau _i \leq b}
\parallel B_i \parallel \exp
\left[ \sum _{b < \tau _i \leq t}
\ln \parallel B_i \parallel \right]
\exp [ - \ln M(I+
\varepsilon ) t ] \leq  $$
$$\leq \prod _{0 < \tau _i \leq b}
\parallel B_i \parallel
\exp [i(t,b) \ln M ]
\exp [ - \ln M(I+
\varepsilon ) t ] \leq  $$
$$\leq \prod _{0 < \tau _i \leq b}
\parallel B_i \parallel
\exp [\ln M  (I +\varepsilon)(t-b)]
\exp [ - \ln M(I+
\varepsilon ) t ]  =  $$
$$= \prod _{0 < \tau _i \leq b}
\parallel B_i \parallel
\exp [- \ln M  (I + \varepsilon )b]
\leq
\prod _{0 < \tau _i \leq b}
\parallel B_i \parallel .
 $$
Thus $\Psi _0 (t) $ is in ${\bf L}_{\infty}$
since
$$ \parallel \Psi _0 (t) \parallel
\leq \prod _{0 < \tau _i \leq b}
( \parallel B_i \parallel  + 1 )  $$
for $t \leq b.$
As
$$ \dot{\Psi} _0 (t) = - \ln M (I + \varepsilon)
\Psi _0 (t) $$
for almost all $t$,
then the derivative of $\Psi _0 $ is also essentially bounded on $[0, \infty)$.

By construction $\Psi_0$ also satisfies the impulsive conditions
$$\Psi _0 (\tau _i) = B_i \Psi _0 (\tau _i - 0).$$
Thus the columns of $\Psi _0$ belong to ${\bf D}_{\infty}.$
Hence the columns of ${\cal M} \Psi _0 $ belong to
${\bf L}_{\infty}$.

If $Y(t,0)$ is the solution of the homogeneous
impulsive equation ${\cal M} y = 0 $, (2), $y(0) = E_n$, then
$ z(t) = \Psi _0 (t) - Y(t,0) $
is the solution of the problem
$$ {\cal M} z = {\cal M} \Psi _0,~
z(\tau _i ) = B_i (\tau _i - 0),~ z(0) = 0.$$
Therefore by Lemma 3.3
$$ Y(t,0) = \Psi _0 (t) -
\int_0^t {C_{{\cal M}} (t,s) ({\cal M} \Psi _0)(s) ds}.$$
The operator $C_{{\cal M}}$ continuously acts from
${\bf L}_{\infty} $ into ${\bf D}_{\infty}$
and the columns of ${\cal M} \Psi _0 $ belong to
${\bf L}_{\infty}$.
 Hence the columns of $Y(t,0)$
belong to ${\bf D}_{\infty},$ therefore $Y(t,0)$
is bounded on the semi-axis $[0, \infty)$
$$ \parallel Y(t,0) \parallel \leq N_0 ,~ t \geq 0.$$
Thus the equality $X(t,0) = \exp {(- \nu t)} Y(t,0) $
implies
$$ \parallel X(t,0) \parallel \leq N_0
 \exp {(- \nu t)}
 ,~ t \geq 0.$$

Now we shall prove that the same estimate is valid
for the Cauchy matrix X(t,s).
The matrix $X(t,s)$ is the solution of the "$s$ -
curtailed" equation (6),(7),(2), therefore by repeating
the proof for this equation one easily obtains
$$ \parallel X(t,s) \parallel \leq N_s
 \exp {[- \nu _s   (t-s)]}
 , t \geq s.$$

It remains to show that $N$ and $\nu$ can be chosen
independently of $s$.
The space ${\bf L}_{\infty}$ contains functions vanishing on
$[0,s)$ , therefore
$$ \parallel C \parallel _{{\bf L}_{\infty} ([s, \infty))
\rightarrow {\bf D}_{\infty}([s, \infty)) } \leq P. $$
The constant $Q$ also does not depend on $s$.
Thus
$$ \parallel {\cal T} C \parallel _{{\bf L}_{\infty} ([s, \infty))
\rightarrow {\bf L}_{\infty}([s, \infty)) } \leq PQ
[\exp {(\nu \delta)} - 1] + P \nu = q < 1 $$
for the same constant $\nu$ .
Here $s \geq 0$ is arbitrary and the values $\nu, q$
do not depend on $s$.
If $\nu$ is chosen small enough for $q < 1$
then the inverse operator in ${\bf L}_{\infty}$
$$(E + {\cal T} C ) ^{-1} $$ exists
and
$$ \parallel (E + {\cal T} C )^{-1}
\parallel _{{\bf L}_{\infty} ([s, \infty))
\rightarrow {\bf L}_{\infty}([s, \infty)) } \leq
\frac{q}{1-q}
$$
for any $s \geq 0.$

Hence (see the proof of Lemma 4.2)
$$ \parallel C_{\cal M} \parallel _{{\bf L}_{\infty} ([s, \infty))
\rightarrow {\bf D}_{\infty} ([s, \infty)) } \leq $$ $$ \leq
\parallel C \parallel _{{\bf L}_{\infty} ([s, \infty))
\rightarrow {\bf D}_{\infty}([s, \infty)) }
\parallel (E - {\cal T} C )^{-1}
\parallel _{{\bf L}_{\infty} ([s, \infty))
\rightarrow {\bf L}_{\infty}([s, \infty)) }.
$$
Thus we have obtained
\begin{eqnarray}
\parallel C_{\cal M} \parallel _{{\bf L}_{\infty} ([s, \infty))
\rightarrow {\bf D}_{\infty}([s, \infty)) } \leq
\frac{qP}{1-q}
\end{eqnarray}
for any $s \geq 0.$

We introduce the functions $\Psi _s (t), ~ s \geq 0, ~ t \geq s $
$$
\Psi _s (t) =
 \left\{ \begin{array}{ll} E_n ,
& \mbox{if~} \tau _i \not\in
(s,t]  ,~ i=1,2, \dots \\
\prod _{s <  \tau _j \leq t} {B_j}, &
\mbox{if~} M \leq 1, \exists \tau _i \in (s,t] , \\
\exp \{ - \ln {M} (I + \varepsilon ) (t-s) \}
\prod _{s <  \tau _j \leq t} {B_j}, &
\mbox{otherwise}.
\end{array}
\right.  $$
As in the case $s=0$
one easily obtains
$$ \parallel \Psi _s (t) \parallel \leq
\prod _{s < \tau _i \leq b }
(\parallel B_i \parallel + 1) \leq
$$ $$ \leq
\prod _{0 < \tau _i \leq b }
(\parallel B_i \parallel + 1)  = U
$$
for every $t \geq s \geq 0$.
Hence $\Psi _s \in {\bf L}_{\infty} . $
Similarly  one obtains that the derivative of $\Psi _s (t) $
in $t$ is  in ${\bf L}_{\infty}$.
Evidently $\Psi _s $ satisfies the impulsive
conditions (2).
Therefore the columns of
$\Psi _s $ are in ${\bf D}_{\infty}$.

By Lemma 3.3 the Cauchy matrix $Y(t,s)$
of the impulsive equation $${\cal M} y = f, ~
y(\tau _i) = B_i y(\tau _i - 0)$$
can be presented as
$$ Y(t,s) = \Psi _s (t) -
\int_s^t {C_{{\cal M}} (t,s) ({\cal M} \Psi _s)(s) ds}.$$
Therefore the inequality (17) implies
$$ \parallel Y(t,s) \parallel
\leq U + q P U \parallel {\cal M} \parallel _{{\bf D}_{\infty}
 \rightarrow {\bf L}_{\infty}}/ (1-q)$$
for any $t,s \geq 0.$
Choosing
$$ N =  q P U \parallel {\cal M} \parallel _{{\bf D}_{\infty}
 \rightarrow {\bf L}_{\infty}} / (1-q) + U $$
we obtain the required estimate
$$ \parallel X(t,s) \parallel
\leq N exp (- \nu (t - s)).$$
The proof of the theorem is complete.

$$ $$
The following example shows that the bounded delay condition
$t - h(t) < \delta $ in Theorem 4.1 is essential.
$$ $$

\underline{\sl Example.}
The solution of the scalar equation
$$ \dot{x} (t) + x(t) - x(0) = f(t) $$
is bounded for any bounded $f$ but $x(t) \equiv 1$
is the solution of the corresponding homogeneous equation.
$$ $$

\section{Exponential stability}
\begin{guess}
 Suppose the hypotheses of Theorem
4.1 hold.
 Then the impulsive equation (1),(2),(3) is exponentially
stable .
\end{guess}

{\sl Proof.}
We immediately obtain the exponential stability
with respect to the initial value from Theorems 3.1 and 4.1
since
$$ \parallel x(t) \parallel =
\parallel X(t,0) x(0) \parallel
\leq N exp(- \nu t) \parallel x(0) \parallel $$
for $r \equiv \varphi \equiv 0.$

As $t - h_i (t) \leq \delta $ then
$$ \varphi [h_i (t)] = 0, ~ t> \delta, ~ i = 1, \dots, m. $$
Therefore for $r \equiv 0, x(0) = 0 $
$$ \parallel x(t) \parallel =
\parallel \int_0^t {X(t,s) \sum _{i=1}^m
{A_i (s) \varphi [h_i(s)] } ds} \parallel = $$
$$ = \parallel \int_0^{\delta} {X(t,s) \sum _{i=1}^m
{A_i (s) \varphi [h_i(s)] } ds} \parallel \leq $$
$$ \leq \int_0^{\delta} { \parallel X(t,s) \parallel \sum _{i=1}^m
{\parallel A_i (s) \parallel
\parallel \varphi [h_i(s)] \parallel  } ds} \leq $$
$$ \leq \sup _{t \geq 0}
\sum _{i=1}^m
{\parallel A_i (t) \parallel}
\cdot \sup _{s<0} \parallel \varphi (s) \parallel
\cdot \int_0^{\delta} {\exp [- \nu (t-s)] ds} \leq $$
$$ \leq \sup _{t \geq 0}
\sum _{i=1}^m
{\parallel A_i (t) \parallel}
\cdot \sup _{s<0} \parallel \varphi (s) \parallel
\cdot \exp (- \nu t) [\exp (\nu \delta ) + 1] .  $$

Denoting
$$  N_0 = \sup _{t \geq 0}
\sum _{i=1}^m
{\parallel A_i (t) \parallel}
\cdot [\exp (\nu \delta ) + 1]   $$
 we obtain
$$ \parallel x(t) \parallel \leq N_0 \exp (- \nu t)
\cdot \sup _{s<0} \parallel \varphi (s) \parallel , $$
i.e. the equation (1),(2)(3) is exponentially stable
with respect to the initial function. The proof of the theorem
is complete.

$$ $$

Theorem 5.1 can be used for obtaining stability results.
For instance the following assertion is valid.
\begin{guess}
Suppose  the hypotheses (a1) -- (a6) hold and
there exist positive constants
$\delta, \zeta , \rho $ such that $t - h_k (t) < \delta,
{}~\zeta  \leq \tau _{i+1} - \tau _i \leq \rho, ~i = 1,2 , \dots,
$
$$\gamma = \sup _i \parallel B_i \parallel < 1 , $$
and
\begin{eqnarray}
\sum _{k=1}^m vraisup _{t>0}
\parallel A_k (t) \parallel \left[  \frac{1}{\alpha} \exp (- \alpha \rho)
+ \rho \right] < 1,
\end{eqnarray}
where
\begin{eqnarray}
\alpha = - \frac{1}{\zeta} \ln {\gamma} .
\end{eqnarray}
Then all solutions of the impulsive equation (1),(2)
are exponentially stable.
\end{guess}

For proving this result we consider an auxiliary ordinary
impulsive equation
\begin{eqnarray}
\dot{x} (t) = f(t),~ t \geq 0, ~ x (\tau _i) = B_i x(\tau _i - 0),
\end{eqnarray}
where $\lim _{i \rightarrow \infty} \tau _i = \infty $.
Let $C_0 (t,s) $ be a Cauchy matrix of (20).
We need the following result for $C_0$.

\begin{uess}
Suppose there exist positive constants
$\rho$ and $\zeta$ such that $$\zeta \leq \tau _{i+1} - \tau _i \leq \rho .$$
Then
$$ \parallel C_0 (t,s) \parallel \leq
\left\{ \begin{array}{ll}
\exp [- \alpha (t-s)], &  \mbox{if~~} t-s > \rho, \\
1, & \mbox{otherwise,}
\end{array}
\right. $$
and for any $t>p, z \in {\bf L}_{\infty}$
the following inequality holds
\begin{eqnarray}
\parallel \int_0^t C_0(t,s) z(s)ds \parallel _{{\bf L}_{\infty}}
\leq \left[ \frac{1}{\alpha} \exp ( - \alpha \rho ) + \rho \right]
\parallel z \parallel _{{\bf L}_{\infty}} .
\end{eqnarray}
Here $\alpha$ is defined by (19).
\end{uess}

{\sl Proof.}
The Cauchy matrix $ C_0 (t,s)$ is the solution of the problem
$$ \dot{x} = 0,~ x(\tau _i) =B_i x(\tau _i - 0), ~ x(s) = E_n,
t \in [s, \infty ). $$
Therefore
$$C_0 (t,s) = \left\{  \begin{array}{ll}
\prod _{s< \tau _i \leq t} B_i, &\mbox{if there exists~}
\tau _i \in (s,t] , \\
E_n, & \mbox{otherwise.}
\end{array}
\right. $$

Let $t -s > \rho$.
Since in the case there exists $i$  such that $s < \tau _i \leq t $
then
$$\parallel C_0(t,s) \parallel \leq \exp \left[
\sum_{s < \tau_i \leq t } \ln \parallel B_i \parallel \right] \leq $$
$$ \leq \exp \left[ (t-s) \frac{i(t,s)}{t-s} \ln(\sup _i \parallel B_i
\parallel ) \right] \leq
\exp \left[\frac{1}{\zeta} (t-s)
\ln {\gamma} \right] = \exp[- \alpha (t-s)]. $$

If $t-s < \rho$ and there exists $\tau _i \in (s,t],$
then the above estimate of $\parallel~C_0(t,s)~\parallel $
is also valid.
As $\ln \gamma < 0$, then
$\parallel C_0 (t,s) \parallel \leq  1.$

If there is no $\tau _i$ belonging to $(s,t]
$ , then $\parallel C_0 (t,s) \parallel = 1.$
Thus we have obtained the estimate of $C_0 (t,s).$

Hence for each $z \in {\bf L}_{\infty}$
$$ \parallel \int_0^t {C_0 (t,s) z(s) ds } \parallel  \leq
\int_0^{t - \rho} {\parallel C_0 (t,s) \parallel
\parallel z(s) \parallel ds } + $$
$$+ \int_{t - \rho}^t {\parallel C_0 (t,s) \parallel
\parallel z(s) \parallel ds }  \leq
\left\{ \int _0^{t - \rho} \exp [- \alpha (t-s)] ds
 + \int_{t -\rho}^t ds \right\}
\parallel z \parallel _{{\bf L}_{\infty}}
= $$
$$ = \left[ \frac{1}{\alpha} \exp (- \alpha \rho) - \frac{1}{\alpha}
\exp (- \alpha t) + \rho \right]
\parallel z \parallel _{{\bf L}_{\infty}} \leq
\left[ \frac{1}{\alpha} \exp (- \alpha \rho)
+ \rho \right]
\parallel z \parallel _{{\bf L}_{\infty}} , $$
which completes the proof of the lemma.

$$ $$

\underline{\sl Proof of Theorem 5.2.}
First we shall make the remark that will be used below.
For the exponential estimation of the Cauchy matrix one
can assume $A_k (t) \equiv 0$ for $t < \rho$.
In fact by Lemma 3.3 if coefficients of two equations
(1),(2) distinguish only on the finite segment [0,b]
then their Cauchy matrices coincide for $s>b$.
Therefore either both of these matrices or none of them
have an exponential estimate.

The proof is based on Theorem 5.1.
Precisely, we shall establish that for any
$f \in {\bf L}_{\infty}$ all solutions of the problem
$$ \dot{x} (t) + \sum _{k=1}^m A_k (t) x[h_k(t)]= f(t), $$
$$x(\tau_i) = B_i x(\tau _i - 0) $$
$$x(0) = 0, $$
$$x(\xi) = 0, \mbox{~if~} ~ \xi < 0, $$
are in ${\bf D}_{\infty}.$

Let $C_0 (t,s)$ be the Cauchy matrix of the equation (20).
Then by substituting
\begin{eqnarray}
x(t) = \int_0^t {C_0 (t,s) z(s) ds } , ~ \dot{x} = z
\end{eqnarray}
we obtain the equation
\begin{eqnarray}
z(t) + \sum_{k=1}^m { \int_0^{h_k^+ (t)}
{A_k (t) C_0 (h_k(t), s) z(s) ds }} = f(t).
\end{eqnarray}

We shall prove that the equation (23) has a solution in
${\bf L}_{\infty}$
for any
$f \in {\bf L}_{\infty}.$
To this end we consider the operator
$$ (Hz)(t) =
\sum_{k=1}^m { \int_0^{h_k^+ (t)}
{A_k (t) C_0 (h_k(t), s) z(s) ds }} $$
and estimate its norm in ${\bf L}_{\infty}$.

Lemma 5.1 and the remark made in the beginning of the proof
imply
$$\parallel H \parallel _{{\bf L}_{\infty} \rightarrow
 {\bf L}_{\infty}} \leq
\sum _{k=1}^m { vraisup _{t>0}
\parallel A_k (t) \parallel \left[\frac{1}{\alpha}
 \exp (- \alpha \rho ) + \rho \right]} .$$
By the inequality (18) $
\parallel H \parallel _{{\bf L}_{\infty} \rightarrow
 {\bf L}_{\infty}} < 1. $
Hence for any $f \in {\bf L}_{\infty}$
the solution $z$ of the equation (23) is in
${\bf L}_{\infty}$.
By (22) and (21) we obtain that
$x \in {\bf L}_{\infty}$ and
$\dot{x} = z \in {\bf L}_{\infty}$.
Thus
$x \in {\bf D}_{\infty}$.

Applying Theorem 5.1 completes the proof of the theorem.
$$ $$

\underline{\sl Remark.} One can easily see that
under the hypotheses of Theorem 5.2 the corresponding equation
without impulses may be unstable.
Thus Theorem 5.2 can be treated as a
stabilization scheme (see also [8]).
$$ $$

\underline{\sl Example.}
The scalar impulsive equation
$$ \dot{x} - a x(t-h) = f,~ a > 0 , $$
$$x(i) = b x(i-1), i= 1,2, \dots , $$
is exponentially stable if
$\mid b \mid < 1 $ and
$$ a \left(1 - \frac{\mid b \mid}{
 \ln {\mid b \mid}} \right) < 1. $$
At the same time the corresponding delay
differential equation
$$ \dot{x} - a x (t-h) = f $$
 is not stable.

\end{document}